\documentclass[pdflatex,sn-mathphys-num]{sn-jnl}

\usepackage{graphicx}%
\usepackage{multirow}%
\usepackage{amsmath,amssymb,amsfonts}%
\usepackage{amsthm}%
\usepackage{mathrsfs}%
\usepackage[title]{appendix}%
\usepackage{xcolor}%
\usepackage{textcomp}%
\usepackage{manyfoot}%
\usepackage{booktabs}%
\usepackage{comment}
\usepackage{algorithm}%
\usepackage{algorithmicx}%
\usepackage{algpseudocode}%
\usepackage{listings}%
\usepackage{bm}

\theoremstyle{thmstyleone}%
\newtheorem{theorem}{Theorem}[section]
\newtheorem{lemma}[theorem]{Lemma}

\theoremstyle{thmstyletwo}%
\newtheorem{remark}{Remark}%
\newtheorem{assumption}{Assumption}
\newtheorem{notation}{Notation}

\theoremstyle{thmstylethree}%
\newtheorem{definition}{Definition}%

\begin{document}

\title[Iterative Bayesian Robbins--Monro Sequence]{An Iterative Bayesian Robbins--Monro Sequence}


\author[1]{\fnm{Siwei} \sur{Liu}}\email{sl2049@cam.ac.uk}

\author[1]{\fnm{Ke} \sur{Ma}}\email{km834@cam.ac.uk}

\author*[1]{\fnm{Stephan M.} \sur{Goetz}}\email{smg84@cam.ac.uk}

\affil*[1]{\orgdiv{Department of Engineering}, \orgname{University of Cambridge}, \orgaddress{\street{Trumpington Street}, \city{Cambridge}, \postcode{CB2 1TN}, \country{United Kingdom}}}


\abstract{This study introduces an iterative Bayesian Robbins--Monro (IBRM) sequence, which unites the classical Robbins--Monro sequence with statistical estimation for faster root-finding under noisy observations. Although the standard Robbins--Monro method iteratively approaches solutions, its convergence speed is limited by noisy measurements and naivety to any prior information about the objective function. The proposed Bayesian sequence dynamically updates the prior distribution with newly obtained observations to accelerate convergence rates and robustness. The paper demonstrates almost sure convergence of the sequence and analyses its convergence rates for both one-dimensional and multi-dimensional problems. We evaluate the method in a practical application that suffers from large variability and allows only a few function evaluations, specifically estimating thresholds in noninvasive brain stimulation, where the method is more robust and accurate than conventional alternatives. Simulations involving 25,000 virtual subjects illustrate  reduced error margins and decreased outlier frequency with direct impact on clinical use.}



\maketitle

\section{Introduction}

Finding the roots of functions is a frequent problem across science and engineering, such as determining physical equilibria or optimising statistical models \cite{bayin2018mathematical,mussardo2010statistical}. Traditional analytical approaches assume that functions are well-defined and analytically tractable. However, real-world problems frequently involve unknown, complex, or noisy functions. Analytical methods are impractical and deterministic methods are typically insufficient in such cases.
The Robbins--Monro sequence has emerged to solve such challenges \cite{robbins1951stochastic}, specifically designed for functions that can only be observed through noisy measurements. 

Assume we are interested in a nondecreasing continuous function $f(x)$ which can only be observed with a bounded zero-mean noise $\varepsilon$, that is, we can only observe 
$y_i$ with $y_i=f(x_i)+\varepsilon_i$ \cite{ruppert1985newton}. Our objective is to find a specific target value $x_{\textrm{t}}$ such that $f(x_{\textrm{t}})=y_{\textrm{t}}$. It is usually assumed that the equation $f(x_{\textrm{t}})=y_{\textrm{t}}$ has a unique root $x_{\textrm{t}}$ in the region of interest \cite{blum1954approximation,ruppert1988efficient}.

Then, starting with an initial guess $x_1$, the algorithm generates a sequence of estimates $\{x_i\}_{i\geq 1}$ according to 
\[
x_{i+1} = x_i - s_i (y_i-y_{\textrm{t}}),
\]
where $\{s_i\}_{i\geq 1}$ is a sequence of positive step sizes \cite{robbins1951stochastic}. To ensure that the Robbins--Monro sequence converges to the true root $x_{\textrm{t}}$, the step size sequence $\{s_i\}_{i\geq 1}$ must satisfy two fundamental conditions \cite{robbins1951stochastic}.

\begin{enumerate}
\item The sum of the whole step-size gain sequence must be infinite:
\[
\sum_{i=1}^\infty s_i = \infty
\]
This condition is essential to guarantee that the sequence takes sufficiently large steps over the course of its iterations to eventually reach the root, even if the initial estimate $x_1$ is far from the true value $x_{\textrm{t}}$ \cite{brent1972optimal,lord1971robbins}. 

\item The sum of the squares of the whole step-size gain sequence must be finite:
\[
\sum_{i=1}^\infty s_i^2 < \infty
\]
This condition is crucial for ensuring that the steps eventually diminish enough so that the sequence converges and does not keep oscillating around the root indefinitely due to the inherent noise in the measuring $f(x_i)$. The finiteness of the sum of squares implies that the variance of the accumulated noise in the estimate remains bounded, which allows the sequence to stabilise and converge \cite{kelley1999iterative,xu2012new}.
\end{enumerate}
One of the most frequently used families of step-size gain sequences is of the form
\[
s_i = \frac{s}{i^\gamma},
\]
where $s$ is a positive constant and $\gamma$ is a parameter in the range $(0.5,1]$ \cite{chen2002robbins}.

If the function $f(x)$ is twice continuously differentiable and strongly monotone in the vicinity of the root, and if the step size sequence is chosen appropriately (e.g., $s_i \propto 1/i$), the Robbins--Monro sequence can achieve an asymptotically optimal convergence rate with respect to the objective function \cite{spall2005introduction,nemirovski2009robust}. 
However, in many practical scenarios, the convergence rate of the standard Robbins--Monro sequence can be relatively slow, e.g., if the step sizes need to decrease slowly to satisfy the convergence conditions in the presence of significant noise or if the starting level of the step sequence is already small so that the one-sided approach towards the root is almost certain but slow \cite{toulis2021proximal}. 

In order to accelerate the convergence rate of the Robbins--Monro sequence, prior information can be introduced into the Robbins--Monro sequence so that the method uses statistical information about where the root might be \cite{liu2024robbins}.
The best estimate of $x_{i+1}$ with the prior known  distribution of $x_\textrm{t}$ ($P_{x_\textrm{t}}(x)$) is
\begin{align}
x_{i+1}=\operatorname{argmax}_x \!\left(P_{x_\textrm{t}}(x) \cdot \mathcal{N}\Bigl(x\, \Big|\, x_i-s_i(y_i-y_\textrm{t}) , c_i^2\Bigr)\right).\label{def}
\end{align} 
The inclusion of prior information accelerates the convergence speed of the Robbins--Monro sequence, especially in the early stage (when there is an accurate prior). However, when the accuracy of the prior distribution is not that high, the convergence speed (w.r.t.\! the number of iterations) will be slow. In order to solve this problem, we propose a new version of the Robbins--Monro sequence with Bayesian learning. 

The structure of this paper is the following: Section \ref{sec2} presents the formal development of our new sequence, the iterative Bayesian Robbins--Monro (IBRM) sequence, and its theoretical foundations. Section \ref{onedim} discusses the convergence and analyses the convergence rate of the one-dimensional new Bayesian Robbins--Monro Sequence. Section \ref{multi} extends the sequence and analysis to multi-dimensional scenarios. Section \ref{tms} applies the new sequence to the estimation of motor thresholds in transcranial magnetic stimulation (TMS) treatment, which needs to measure a response level with few measurements despite high trial-to-trial measurement variability \cite{wang2023three}. Finally, Section \ref{con} concludes the paper.

\section{Sequence development} \label{sec2}
\begin{notation}
$\mathcal{N}(x \,|\, \alpha, \beta^2)$ is the normal distribution with $\alpha$ mean and $\beta^2$ variance.\\ 
\end{notation}

Following Bayes's rule, we can have
\begin{align}
P(x_\textrm{t} \mid x_{i+1}) \propto P_{x_\textrm{t}}(x)  \cdot P(x_{i+1} \mid x_i, s_i, \mathrm {Model}), \quad \forall i \geq 1.
\end{align}
Now, we assume the prior distribution $P_{x_\textrm{t}}(x)$ can be updated with the current information. That is, when we have some information about $x_i$, we can use $P(x_\textrm{t} \mid x_{i})$ to be the prior distribution and update the sequence. In this case, we set the step size to be constantly $s$. Then, the sequence becomes 
\begin{align}
P(x_\textrm{t} \mid x_{i+1}) \propto P(x_\textrm{t} \mid x_{i})  \cdot P(x_{i+1} \mid x_i, s, \mathrm {Model}), \quad \forall i \geq 1.
\end{align}
If we start the investigation from the first step, we can get
\begin{align}
P(x_\textrm{t} \mid x_{2}) &\propto P_{x_\textrm{t}}(x) \cdot P(x_{2} \mid x_1, s, \mathrm {Model}) \notag \\
&\propto P_{x_\textrm{t}}(x) \cdot \mathcal{N}\Bigl(x\, \Big|\, x_1-s(y_1-y_\textrm{t}) , c_1^2\Bigr). \label{eqpr1}
\end{align}
In this way, Eq.\ (\ref{eqpr1}) can be the updated (unnormalised) prior on the second iteration. We repeat this process iteratively and get
$$
P(x_\textrm{t} \mid x_{i}) 
\propto P_{x_\textrm{t}}(x) \cdot \prod_{k=1}^{i-1} \mathcal{N}\Bigl(x\, \Big|\, x_k-s(y_k-y_\textrm{t}) , c_k^2\Bigr).
$$
and 
$$
x_{i+1}= \operatorname{argmax}_x \left(P_{x_\textrm{t}}(x) \cdot \prod_{k=1}^i \mathcal{N}\Bigl(x\, \Big|\, x_k-s(y_k-y_\textrm{t}) , c_k^2\Bigr)\right)
$$
with $c_k=c/k$.
We call this new sequence the iterative Bayesian Robbins--Monro sequence (IBRM).
In the high-dimensional case, this sequence becomes 
$$
\bm{x}_{i+1}= \operatorname{argmax}_{\bm{x}} \left(P_{\bm{x}_\textrm{t}}(\bm{x}) \cdot  \prod_{k=1}^i \mathcal{\bm{N}}\Bigl(\bm{x}\, \Big|\, \bm{x}_k-s(\bm{y}_k-\bm{y}_\textrm{t}) , \bm{C}_k \Bigr)\right),
$$ 
where $\bm{C}_k$ is the covariance matrix.

\begin{notation}
In the following analysis, we call $P_{\bm{x}_\textrm{t}}(\bm{x})$ as well as $P_{x_\textrm{t}}(x)$ prior distribution and $\mathcal{N}(x\, |\, x_k-s(y_k-y_\textrm{t}) , c_k^2)$ as well as $\mathcal{\bm{N}}(\bm{x}\, |\, \bm{x}_k-s(\bm{y}_k-\bm{y}_\textrm{t}) , \bm{C}_k )$ RM distribution.
\end{notation}

\section{One-dimensional iterative Bayesian Robbins--Monro sequence}\label{onedim}

\subsection{Convergence analysis} 

Firstly, we investigate the case where the prior distribution is a normal distribution.

\begin{notation} \label{n2} 
In all the following proof, $\prod_{u=m}^{n} Q(u) =1$ if $m>n$ for any arbitrary function $Q$.
\end{notation}

\begin{notation}
All norms we use in this paper are 2-norms, that is, $\| \cdot\|=\| \cdot\|_2$.
\end{notation}

\begin{assumption} \label{a1}
Assume w.l.o.g.\ that the variance of the RM distribution follows $1/k^2$ in the following proof, that is, $c_k=1/k$ and $c=1$.
\end{assumption}
\begin{assumption}
The noise variables $\{\varepsilon_i\}$ are independently drawn from a distribution that is symmetric about zero, i.e., the probability density satisfies $P_{\varepsilon_i}(x) = P_{\varepsilon_i}(-x)$ for all $x \in \mathbb{R}$ and $i \in \mathbb{Z}^+$.
\end{assumption}
\begin{assumption} 
The target point $x_\textrm{t}$ is finite.
\end{assumption}
\begin{assumption} \label{a6}
The underlying function $f$ of the Robbins--Monro sequence is an increasing function.
\end{assumption}
\textbf{In the absence of further clarification, we assume by default that Assumptions \ref{a1}--\ref{a6} hold for all the proofs in Section \ref{onedim}.}
\\
\\
\begin{theorem}\label{thm1}
When the prior distribution is normal, the iterative Bayesian Robbins--Monro sequence generated by 
$$
x_{i+1}= \operatorname{argmax}_x \left(\mathcal{N} ( x \mid \mu_n, \sigma_n^2) \cdot  \prod_{k=1}^i \mathcal{N}\Bigl(x\, \Big|\, x_k-s(y_k-y_\textrm{t}) ,\frac{1}{k^2} \Bigr)\right),
$$ 
where $ \mu_n, \sigma_n$, and $s$ are finite constants with $ \sigma_n >0$, is convergent a.s.\! (almost surely).
\begin{proof}
We want to prove that the iterative Bayesian Robbins--Monro sequence is indeed a Robbins--Monro sequence.

We know that 
\begin{align}
x_{i}&= \operatorname{argmax}_x \left(\mathcal{N} ( x \mid \mu_n, \sigma_n^2) \cdot \prod_{k=1}^{i-1} \mathcal{N}\Bigl(x\, \Big|\, x_k-s(y_k-y_\textrm{t}) ,\frac{1}{k^2} \Bigr)\right).\\
x_{i+1}&= \operatorname{argmax}_x \left(\mathcal{N} ( x \mid \mu_n, \sigma_n^2) \cdot \prod_{k=1}^{i-1} \mathcal{N}\Bigl(x\, \Big|\, x_k-s(y_k-y_\textrm{t}) ,\frac{1}{k^2} \Bigr) \cdot \mathcal{N}\Bigl(x\, \Big|\, x_i-s(y_i-y_\textrm{t}) ,\frac{1}{i^2} \Bigr)\right) \label{eqi}
\end{align}
and the product of probability density functions (pdf) of normal distributions is the unnormalised pdf of another normal distribution. Therefore, we can see that $x_{i+1}$ lies between $x_i$ and $x_i-s(y_i-y_\textrm{t})$.

Then, it is sufficient to find the equivalent step size sequence $\{s_i, \ i=1,2,...\}$ such that $x_{i+1}=x_i-s_i(y_i-y_\textrm{t})$ and prove that this step size sequence satisfies
$$
\begin{aligned}
\sum_{i=1}^\infty s_i &=\infty,\\
\sum_{i=1}^\infty s_i^2 &< \infty
\end{aligned}
$$
to show that IBRM is indeed a standard Robbins--Monro sequence and hence is convergent a.s. 

The variance of
$\prod_{k=1}^{i-1} \mathcal{N}\Bigl(x\, \Big|\, x_k-s(y_k-y_\textrm{t}) ,1/k^2 \Bigr)$
is 
$$
\begin{aligned}
\frac{1}{\sum_{k=1}^{i-1} k^2 }&= \frac{1}{\mathcal{O}(1) \cdot \int_{k=0}^{i-1} k^2 \,\textrm{d}k} = \mathcal{O}(1) \cdot \frac{1}{i^3}. 
\end{aligned}
$$
Since $ \mu_n $ and $ \sigma_n$ are finite constants, the variance of 
$$
\mathcal{N} ( x \mid \mu_n, \sigma_n^2) \cdot \prod_{k=1}^{i-1} \mathcal{N}\Bigl(x\, \Big|\, x_k-s(y_k-y_\textrm{t}) ,\frac{1}{k^2} \Bigr)
$$
is also $\mathcal{O}(1) \cdot 1/i^3$.
Therefore, when looking at Eq.\! (\ref{eqi}), we can know that 
$$
x_{i+1} =  \frac{\mathcal{O}(1) \cdot i^3  x_i + i^2(x_i-s(y_i-y_\textrm{t}))}{\mathcal{O}(1) \cdot i^3+ i^2}.
$$
Therefore, 
$$
s_i= \frac{ i^2 s}{\mathcal{O}(1) \cdot i^3+ i^2} =\mathcal{O}(1) \cdot \frac{ 1}{i},
$$
which satisfies the step size conditions required by the standard Robbins--Monro sequence.

Therefore, the iterative Bayesian Robbins--Monro sequence is indeed a Robbins--Monro sequence and hence is a.s.\ convergent.
\end{proof}
\end{theorem}

\begin{definition}\label{def1}
Let $x_1^{\mathrm{r}}, x_2^{\mathrm{r}},...$ be the standard Robbins--Monro sequence and $s_1, s_2,...$ be step sizes for the standard Robbins--Monro sequence. Then, the standard Robbins--Monro sequence follows 
$$
x_{i+1}^{\mathrm{r}}=x_i^{\mathrm{r}}-s_i(y_i^{\mathrm{r}}-y_\textrm{t}).
$$

In the standard Robbins--Monro sequence, the target point we want to find is also $x_\textrm{t}$. Assume $f(x_\textrm{t})=y_\textrm{t}$, and $\sum_{i=1}^\infty s_i=\infty$, $\sum_{i=1}^\infty s_i^2 < \infty$. Moreover, $y_i^{\mathrm{r}}=y(x_i^{\mathrm{r}})$ is the noisy measurement of $f(x)$ at $x_i^{\mathrm{r}}$ with $y_i^{\mathrm{r}}= f(x_i^{\mathrm{r}})+\varepsilon_{i}^{\mathrm{r}}$ and $\{\varepsilon_{i}^{\mathrm{r}}\}_{i \geq 1}$ are bounded independent random variables with 0 mean as well as $d^2$\! variance ($d<\infty$). \\
\end{definition}

\begin{lemma}\label{lem:important}
When $f(x)=bx$, i.e., linear, ($b>0$), for any finite $h \in \mathbb{Z}^{+},$
$$
\begin{aligned}
\prod_{k=h}^\infty (1-b s_k) = 0 \quad a.s.
\end{aligned}
$$
\begin{proof}
\begin{align}
x_{i+1}^{\mathrm{r}} &=x_i^{\mathrm{r}} -s_i (b(x_i^{\mathrm{r}}-x_\textrm{t})+\varepsilon_{i}^{\mathrm{r}})\notag\\
x_{i+1}^{\mathrm{r}}-x_\textrm{t} &=(1-bs_i)(x_i^{\mathrm{r}}-x_\textrm{t}) +s_i \varepsilon_{i}^{\mathrm{r}}\notag\\
&=\prod_{k=h}^{i}(1-b s_k)(x_h^{\mathrm{r}}-x_\textrm{t})\label{eq1.4} \\
&+\sum_{k=h}^{i}s_k\varepsilon_{k}^{\mathrm{r}}\prod_{u=k+1}^{i} (1-bs_u)\label{eq1.5}
\end{align}
where $h \in \mathbb{Z}^{+}$ is an arbitrary number less than $i$.

Here, we rephrase 
$$\sum_{k=1}^{i-1}s_k\varepsilon_{k}^{\mathrm{r}}\prod_{u=k+1}^{i} (1-bs_u)+ s_i\varepsilon_{i}^{\mathrm{r}}$$ 
into 
$$\sum_{k=1}^{i}s_k\varepsilon_{k}^{\mathrm{r}}\prod_{u=k+1}^{i} (1-bs_u)$$
with Notation \ref{n2} for simplification. We will repeatedly use this abbreviation in the following proofs.

Since the standard Robbins--Monro sequence is convergent a.s.\ and $E(\varepsilon_{k}^{\mathrm{r}})=0$ $\forall k$,  we can get
$$
\begin{aligned}
&\lim_{i \to \infty}E \big[x_{i+1}^{\mathrm{r}}-x_\textrm{t} \big] \\
=&\prod_{k=h}^{\infty}(1-b s_k)(x_h^{\mathrm{r}}-x_\textrm{t}) =0 \quad a.s.
\end{aligned}
$$
for all finite $h \in \mathbb{Z}^{+}$. 

Moreover, $x_h^{\mathrm{r}}-x_\textrm{t} \neq 0$ a.s.\ since the value of $x_h^{\mathrm{r}}$ is influenced by $h-1$ continuous random variables $\varepsilon_{1}^{\mathrm{r}},\varepsilon_{2}^{\mathrm{r}},...,\varepsilon_{h-1}^{\mathrm{r}}$. Therefore,
\begin{align*}
 \prod_{k=h}^{\infty}(1-b s_k) = 0 \quad a.s.
\end{align*}
for all finite $h \in \mathbb{Z}^{+}$.
\end{proof}
\end{lemma}

\begin{lemma}\label{lem2}
When $f(x)=bx$, i.e., linear, and $b>0$, $$\sum_{k=h}^{\infty} s_k^2 d^2 \prod_{u=k+1}^{\infty}(1-b s_u)^2  = 0 \quad a.s.$$ 
for arbitrary $h \in \mathbb{Z}^{+}$. Note that $d^2$\! is the variance of $\{\varepsilon_i\}$.
\begin{proof}
Calculate $\left(\text{Term } (\ref{eq1.4}) + \text{Term }(\ref{eq1.5})\right)^2$ and take the expectation to get
\begin{align}
E[(x_{i+1}^{\mathrm{r}}-x_\textrm{t})^2]&=\prod_{k=h}^i(1-b s_k)^2(x_h^{\mathrm{r}}-x_\textrm{t})^2 \label{eq2}\\
&+\sum_{k=h}^i s_k^2 d^2 \prod_{{u=k+1 }}^i(1-b s_u)^2 \label{eq3}
\end{align}
since $\{\varepsilon_{i}^{\mathrm{r}}\}_{i \geq 1}$  are independent random variables with 0 mean and $d^2$ variance. As the standard RM sequence is a.s.\! convergent, 
\begin{align*}
\lim_{i \to \infty}E[(x_{i+1}^{\mathrm{r}}-x_\textrm{t})^2]  = 0 \quad a.s.
\end{align*}
\\
Since Terms $(\ref{eq2}), (\ref{eq3}) \geq 0$, then Terms $(\ref{eq2}), (\ref{eq3}) \stackrel{i\to \infty}{\longrightarrow} 0$ a.s.{},
$$
\sum_{k=h}^{\infty} s_k^2 d^2 \prod_{u=k+1 }^{\infty}
(1-b s_u)^2 =0 \quad a.s.
$$
\end{proof}
\end{lemma}

\begin{assumption} \label{a3}
If $f(x)$ is a nonlinear function, then assume $f(x)$ is continuously differentiable with $f'(x)$ being positive and bounded per $M \geq f^{\prime} \geq m > 0$.
\end{assumption}

\begin{assumption} \label{a9} 
The prior distribution $P_{x_\textrm{t}}(x)$ is a finite differentiable function that satisfies $|P_{x_\textrm{t}}^{\prime}(x)/ P_{x_\textrm{t}} (x)| \leq J,$ $\forall x \in \mathbb{R}$ for some $J>0$. Furthermore, the prior should be bounded per $ 0< P_{x_\textrm{t}} (x) < \infty$. 
\end{assumption}

\begin{remark}
We may see that if $P_{x_\textrm{t}}(x)$ is a Laplace distribution, Logistic distribution, Cauchy distribution, or Student’s t-distribution, the Assumption \ref{a9} will be satisfied for $P_{x_\textrm{t}}(x)$.
\end{remark}

\begin{theorem}\label{thm2}
When Assumptions \ref{a3}--\ref{a9} holds, the iterative Bayesian Robbins--Monro sequence generated by 
$$
x_{i+1}= \operatorname{argmax}_x \left( P_{x_\textrm{t}}(x) \cdot \prod_{k=1}^i \mathcal{N}\Bigl(x\, \Big|\, x_k-s(y_k-y_\textrm{t}) ,\frac{1}{k^2} \Bigr)\right),
$$ 
where $ \mu_n, \sigma_n$, and $s$ are finite constants with $ \sigma_n >0$, is convergent a.s.
\begin{proof}
Denote
$$
\begin{aligned}
x_{i}^{\textrm{c}}&= \operatorname{argmax}_x \left( \prod_{k=1}^{i-1} \mathcal{N}\Bigl(x\, \Big|\, x_k-s(y_k-y_\textrm{t}) ,\frac{1}{k^2} \Bigr)\right) \\
&=  \operatorname{argmax}_x \left(  \mathcal{N}(x\, \Big|\,x_{i}^{\textrm{c}} ,\sigma^2_i )\right)\\
&= \operatorname{argmax}_x \left(  N_i(x)\right)
\end{aligned}
$$
and
$$
\begin{aligned}
x_{i+1}^{\textrm{c}}&= \operatorname{argmax}_x \left( \prod_{k=1}^{i} \mathcal{N}\Bigl(x\, \Big|\, x_k-s(y_k-y_\textrm{t}) ,\frac{1}{k^2} \Bigr)\right) \\
&=  \operatorname{argmax}_x \left(  \mathcal{N}(x\, \Big|\,x_{i+1}^{\textrm{c}} ,\sigma^2_{i+1} )\right)\\
&= \operatorname{argmax}_x \left(  N_{i+1}(x) \right).
\end{aligned}
$$
\\
\\
When $N_i(x) \cdot P_{x_\textrm{t}}(x)$ attains its maximum at $x_{i}$, 
$$
x_{i}=\operatorname{argmax}_x\big(P_{x_\textrm{t}}(x) \cdot N_i(x)\big),
$$
$$
N_i^{\prime}(x_{i}) P_{x_\textrm{t}}(x_{i})+ N_i(x_{i}) P_{x_\textrm{t}}^{\prime}(x_{i})=0,
$$ 

$$
\frac{N_i^{\prime}(x_{i})}{N_i(x_{i})}=-\frac{P_{x_\textrm{t}}^{\prime} (x_{i})}{P_{x_\textrm{t}}(x_{i})}.
$$
\\
Therefore,
$$
x_i-x_{i}^{\textrm{c}}= \sigma_i^2 \cdot \frac{P_{x_\textrm{t}}^{\prime} (x_{i})}{ P_{x_\textrm{t}}(x_{i})}. 
$$
\\
Following the idea in Theorem \ref{thm1}, we know that 
$$
\sigma_i^2 = \mathcal{O}(1) \cdot \frac{1}{i^3}.
$$
Moreover,
$$\left|\frac{P_{x_\textrm{t}}^{\prime}(x_{i})}{P_{x_\textrm{t}}(x_{i})}\right| \leq J.$$
Therefore,
$$
x_i-x_{i}^{\textrm{c}}= \mathcal{O}(1) \cdot \frac{1}{i^3}.
$$
Similarly,
$$
x_{i+1}-x_{i+1}^{\textrm{c}}= \mathcal{O}(1) \cdot \frac{1}{i^3}.
$$
\\
Moreover, from Theorem \ref{thm1}, we know that 
$$
x_{i+1}^{\textrm{c}}=x_{i}^{\textrm{c}}-s_i (y_i-y_\textrm{t})
$$
with 
$$
s_i =\mathcal{O}(1) \cdot \frac{1}{i}.
$$
\\
Hence, 
$$
\begin{aligned}
x_{i+1}=&x_{i}-s_i (y_i-y_\textrm{t})+\mathcal{O}(1) \cdot \frac{1}{i^3}, \\
x_{i+1}=&x_{i}-s_i \left(f'(z_i)(x_i-x_\textrm{t})\right)+ s_i \varepsilon_i +\mathcal{O}(1) \cdot \frac{1}{i^3}, \\
x_{i+1}-x_\textrm{t}=&(1-s_i f'(z_i))(x_i-x_\textrm{t})+ s_i \varepsilon_i+\mathcal{O}(1) \cdot \frac{1}{i^3}. \\
\end{aligned}
$$
\\
By expanding the equation above, we can get
\begin{align}
\lim_{i \to \infty}x_{i+1}-x_\textrm{t}=&\prod_{k=h}^{\infty}(1- s_k f'(z_k) )(x_h-x_\textrm{t}) \label{eq11} \\
+&\sum_{k=h}^{\infty}s_k\varepsilon_{k} \prod_{u=k+1}^{\infty} (1- s_u f'(z_u)) \label{eq12} \\
+&\mathcal{O}(1) \cdot \sum_{k=h}^{\infty} \frac{1}{k^3} \prod_{u=k+1}^{\infty} (1- s_u f'(z_u)), \label{eq13}
\end{align}
where $h$ is a given large number that satisfies $s_k< 1/M,$ $\forall k \geq h$.
\\
\\
From Lemma \ref{lem:important}, we know that $\prod_{k=h}^{\infty}(1- s_k m )=0$ a.s. Hence, Term (\ref{eq11}) $=0$ a.s.\ as $i \to \infty$.
Moreover, from Lemma \ref{lem2} we have known that 
$$
\sum_{k=h}^{\infty} s_k^2 d^2 \prod_{u=k+1 }^{\infty}(1- s_u m)^2 =0 \quad a.s.
$$
Hence, $E((\ref{eq12})^2)=0$ when $i \to \infty$.
\\
\\
Since $ \sum_{k=1}^{\infty} 1/k^3 < \infty$, for any infinitesimal $p >0$, there exists $Z \in \mathbb{N}$ such that 
$$\left| \sum_{k=Z+1}^{\infty} \frac{1}{k^3}\right| < p.$$
W.l.o.g., assume $Z>h$. Therefore, \\
$$
0 \leq \left|\sum_{k=Z+1}^{\infty} \frac{1}{k^3} \prod_{u=k+1}^{\infty} (1- s_u f'(z_u))\right| \leq \left|\sum_{k=Z+1}^{\infty} \frac{1 }{k^3} \right| < p.
$$
\\
Moreover, when $h \leq k \leq Z$, following the same idea in proving in Term (\ref{eq11}) $=0$, we can get
$$
\sum_{k=h}^{Z} \frac{1}{k^3} \prod_{u=k+1}^{\infty} (1- s_u f'(z_u))=0 \quad a.s.
$$
Therefore, we can get $\left|\textrm{Term (\ref{eq13})}\right| < p$ a.s.\! for any infinitesimal $p > 0$ when $i \to \infty$. 
Hence, 
$$
\lim_{i \to \infty}x_{i+1} = x_\textrm{t} \quad a.s. 
$$
and the iterative Bayesian Robbins--Monro sequence is convergent a.s.\! when Assumptions \ref{a3}--\ref{a9} hold.

\end{proof}
\end{theorem}


    

\subsection{Convergence rate analysis}

\begin{remark}
When the prior distribution is normal, the iterative Bayesian Robbins--Monro sequence shares the convergence rate with the standard Robbins--Monro sequence in the long run. The initial convergence rate of the iterative Bayesian Robbins--Monro sequence depends on the accuracy of the prior distribution.
\end{remark}

Then, we will analyse the convergence rate of the case corresponding to Theorem \ref{thm2}.
\\
\\
\begin{lemma} \label{lem35}
When $f(x)=bx$, i.e., linear ($b>0$),  
$\|E((\ref{eq12})^2)\| \gg \|(\ref{eq13})\|^2 $ when $i \to \infty.$
\begin{proof}
When $f'=b$, we can get
$$
\begin{aligned}
&\ln\left(\prod_{u= k+1}^i(1-b s_u)\right) = \sum_{u= k+1}^i\ln(1-b s_u) =\mathcal{O}(1) \cdot \sum_{u= k+1}^i -\frac{b}{u}\\
=& - \mathcal{O}(1) \cdot b \int _{k+1}^i \frac{1}{u} \, du =\mathcal{O}(1) \cdot b \ln\frac{k}{i} = \mathcal{O}(1) \cdot\ln\frac{k^b}{i^b}.
\end{aligned}
$$
Therefore,
$$
\prod_{u= k+1}^i(1-b s_u) =\frac{k^{bv}}{i^{bv}}.
$$
for some constant $v=\mathcal{O}(1)$. Then, we can see that when $i$ is finite, 
$$
\begin{aligned}
\textrm{(\ref{eq13})} = & \mathcal{O}(1) \cdot \sum_{k=h}^i \frac{1}{k^3} \cdot \frac{k^{bv}}{i^{bv}} = \mathcal{O}(1) \cdot\sum_{k=h}^i  \frac{k^{bv-3}}{i^{bv}} \\
= & \mathcal{O}(1) \cdot \frac{1}{i^{bv}} \int_{h}^i k^{bv-3} \, dk 
= \begin{cases}\mathcal{O}(1) \cdot\ln i/i^2 & \text { if } bv=2, \\
\mathcal{O}(1) \cdot 1/i^2 & \text { otherwise, }\end{cases} \\
\end{aligned}
$$
and 
$$
\|\textrm{(\ref{eq13})}\|^2 \sim \begin{cases}\mathcal{O}(1) \cdot (\ln i)^2/i^4 & \text { if } bv=2, \\
\mathcal{O}(1) \cdot 1/i^4 & \text { otherwise. }\end{cases} \\
$$

Similarly, we can see that 
$$
\begin{aligned}
&\| E( (\ref{eq12})^2) \| = \sum_{k=h}^{i}  s_k^2 d^2 \cdot \prod_{u=k+1}^{i} \left( 1-b s_u \right) ^2\\
= & \mathcal{O}(1) \cdot \sum_{k=h}^i  \frac{k^{2 bv-2}}{i^{2 bv}} =\begin{cases}
\mathcal{O}(1) \cdot \ln i/i & \text { if } bv=0.5, \\
 \mathcal{O}(1) \cdot 1/i & \text { otherwise, }
\end{cases} \\
&
\end{aligned}
$$
Note that 
$$
\frac{\ln i}{i} \gg \frac{1}{i} \gg \frac{(\ln i)^2}{i^4} \gg \frac{1}{i^4}
$$
when $i \to \infty$. Therefore, if $f$ is a linear function,
$\|E((\ref{eq12})^2)\| \gg \|(\ref{eq13})\|^2 $ when $i \to \infty.$

\end{proof}
\end{lemma}

\begin{theorem} \label{thm36}
When Assumptions \ref{a3}--\ref{a9} hold, the iterative Bayesian Robbins--Monro sequence generated by 
$$
x_{i+1}= \operatorname{argmax}_x \left( P_{x_\textrm{t}}(x) \cdot \prod_{k=1}^i \mathcal{N}\Bigl(x\, \Big|\, x_k-s(y_k-y_\textrm{t}) ,\frac{1}{k^2} \Bigr)\right),
$$ 
where $ \mu_n, \sigma_n$, and $s$ are finite constants with $ \sigma_n >0$, has the same convergence rate as the standard Robbins--Monro sequence in the long run.

\begin{proof}
W.l.o.g, assume $x_h=x_h^\textrm{r}$. Comparing with Terms (\ref{eq1.4})--(\ref{eq1.5}), we know that 
$$
\left|\prod_{k=h}^{\infty}(1- s_k m )(x_h-x_\textrm{t})\right| \geq \left| (\ref{eq11}) \right| \geq \left|\prod_{k=h}^{\infty}(1- s_k M )(x_h-x_\textrm{t})\right|.
$$
Moreover, 
$$
\left|\sum_{k=h}^\infty s_k^2 d^2 \prod_{{u=k+1 }}^\infty (1-m s_u)^2\right| \geq \left| E((\ref{eq12})^2) \right| \geq \left|\sum_{k=h}^\infty s_k^2 d^2 \prod_{{u=k+1 }}^\infty (1-M s_u)^2\right|.
$$
Then,
$$
\mathcal{O}(1) \cdot \left|\sum_{k=h}^{\infty} \frac{1}{k^3} \prod_{u=k+1}^{\infty} (1- s_u m)\right| \geq \left| (\ref{eq13}) \right| \geq \mathcal{O}(1) \cdot \left|\sum_{k=h}^{\infty} \frac{1}{k^3} \prod_{u=k+1}^{\infty} (1- s_u M)\right|.
$$
Moreover, by Lemma \ref{lem35},
$$
\left|\sum_{k=h}^\infty s_k^2 d^2 \prod_{{u=k+1 }}^\infty (1-M s_u)^2\right| \gg \mathcal{O}(1) \cdot \left|\sum_{k=h}^{\infty} \frac{1}{k^3} \prod_{u=k+1}^{\infty} (1- s_u m)\right|.
$$
Therefore, for an arbitrary $f$ with a bounded and positive gradient,
$\|E((\ref{eq12})^2)\| \gg \|(\ref{eq13})\|^2 $ when $i \to \infty.$ We can ignore Term (\ref{eq13}) when analysing the convergence rate of the iterative Bayesian Robbins--Monro sequence.
Since the standard Robbins--Monro sequences with $f(x)=mx$ and $f(x)=Mx$ have the same convergence rate, we can say that the iterative Bayesian Robbins--Monro sequence has the same convergence rate as the standard Robbins--Monro sequence in the long run when Assumptions \ref{a3}--\ref{a9} hold.
\end{proof}
\end{theorem}

\section{Multi-dimensional iterative Bayesian Robbins--Monro sequence} \label{multi}

\begin{notation} In all the following proofs, $\prod_{k \in \text{desc}(h_1, h_2)} \bm{Q_k} = \bm{Q_{h_2}} \cdot \bm{Q_{h_2-1}} \cdots \bm{Q_{h_1}} $ if $h_2 > h_1$ and  $h_1,$ $h_2 \in \mathbb{Z}^{+}$ for any arbitrary series of matrix $\{\bm{Q_k}\}$. If $h_1 > h_2$, $\prod_{k \in \text{desc}(h_1, h_2)} \bm{Q_k} = \bm{I} $.
\end{notation}

\begin{assumption} \label{b1}
W.l.o.g., assume the covariance matrix of the RM distribution follows $1/k^2 \cdot \bm{I}$ in the following proof.
\end{assumption}
\begin{assumption}
Assume $\bm{f}: \mathbb{R}^q \rightarrow \mathbb{R}^q$ is continuous differentiable. Moreover, $\nabla\bm{ f}$ is a symmetric p.d.\! matrix that satisfies $ m\bm{I} \prec \nabla\bm{ f} \prec M\bm{I}$, where $0 < m < M$. 
\end{assumption}
\begin{assumption}
The noise variables $\{\bm{\varepsilon}_i\}$ are independently drawn from a distribution that is symmetric about zero, i.e., the probability density satisfies $P_{\bm{\varepsilon}_i}(\bm{x}) = P_{\bm{\varepsilon}_i}(-\bm{x})$ for all $\bm{x} \in \mathbb{R}^q$ and $i \in \mathbb{Z}^+$.
\end{assumption}
\begin{assumption} 
The target point $\bm{x}_\textrm{t}$ is finite.
\end{assumption}
\begin{assumption} \label{b6}
The multi-dimensional prior distribution $P_{\bm{x}_\textrm{t}}(\bm{x})$ is a finite differentiable function that satisfies $||\nabla P_{\bm{x}_\textrm{t}}(\bm{x})/ P_{\bm{x}_\textrm{t}} (\bm{x})|| \leq J,$ $\forall \bm{x} \in \mathbb{R}^q$ for some $J>0$. Furthermore, the prior should satisfy $ \bm{0}< P_{\bm{x}_\textrm{t}} (\bm{x}) < \bm{\infty}$. 
\end{assumption}

\textbf{In the absence of further clarification, we assume by default that Assumptions \ref{b1}--\ref{b6} hold for all the proofs in Section \ref{multi}.}
\\
\\
\begin{theorem}
The multi-dimensional iterative Bayesian Robbins--Monro sequence generated by 
$$
\bm{x}_{i+1}= \operatorname{argmax}_{\bm{x}} \left(P_{\bm{x}_\textrm{t}}(\bm{x}) \cdot  \prod_{k=1}^i \mathcal{\bm{N}}\Bigl(\bm{x}\, \Big|\, \bm{x}_k-s(\bm{y}_k-\bm{y}_\textrm{t}) ,\frac{1}{k^2} \bm{I} \Bigr)\right),
$$ 
where $s$ is a finite positive constant, is convergent a.s..
\begin{proof}
As in Theorem \ref{thm2}, we define 
$$
\begin{aligned}
\bm{x}_{i}^{\textrm{c}}&= \operatorname{argmax}_{\bm{x}} \left( \prod_{k=1}^{i-1} \mathcal{\bm{N}}\Bigl(\bm{x}\, \Big|\, \bm{x}_k-s(\bm{y}_k-\bm{y}_\textrm{t}) ,\frac{1}{k^2} \bm{I} \Bigr)\right) \\
&=  \operatorname{argmax}_{\bm{x}} \left(  \mathcal{\bm{N}}(\bm{x} \, \Big|\,\bm{x}_{i}^{\textrm{c}} ,\bm{\Lambda}_i )\right)\\
&= \operatorname{argmax}_{\bm{x}} \left(  \bm{N}_i(\bm{x})\right)
\end{aligned}
$$
and
$$
\begin{aligned}
\bm{x}_{i+1}^{\textrm{c}}&= \operatorname{argmax}_{\bm{x}} \left( \prod_{k=1}^{i} \mathcal{\bm{N}}\Bigl(\bm{x}\, \Big|\, \bm{x}_k-s(\bm{y}_k-\bm{y}_\textrm{t}) ,\frac{1}{k^2} \bm{I} \Bigr)\right) \\
&=  \operatorname{argmax}_{\bm{x}} \left(  \mathcal{\bm{N}}(\bm{x}\, \Big|\,\bm{x}_{i+1}^{\textrm{c}} ,\bm{\Lambda}_{i+1} )\right)\\
&= \operatorname{argmax}_{\bm{x}} \left(  \bm{N}_{i+1}(\bm{x}) \right).
\end{aligned}
$$
\\
\\
Moreover, we can get 
\[
\bm{\Lambda}_i = \left( \sum_{k=1}^{i-1} k^2 \right)^{-1} \bm{I}=\mathcal{O}(1) \cdot \frac{1}{i^3} \bm{I}.
\]
Then, following the same idea as in Theorem \ref{thm1}, we can derive
$$
\bm{x}_{i+1}^{\textrm{c}}=\bm{x}_{i}^{\textrm{c}}-s_i(\bm{y}_i-\bm{y}_\textrm{t})
$$
with 
$$
s_i= \mathcal{O}(1) \cdot \frac{ 1}{i}
$$
since all the covariance matrices above are isotropic.

When $\bm{N}_i(\bm{x}) \cdot P_{\bm{x}_\textrm{t}}(\bm{x})$ attains its maximum at $\bm{x}_{i}$, 
$$
\bm{x}_{i}=\operatorname{argmax}_{\bm{x}} \big(P_{\bm{x}_\textrm{t}}(\bm{x}) \cdot \bm{N}_i(\bm{x})\big),
$$
$$
\nabla \bm{N}_i(\bm{x}_{i}) P_{\bm{x}_\textrm{t}}(\bm{x}_{i})+ \bm{N}_i(\bm{x}_{i}) \nabla P_{\bm{x}_\textrm{t}}(\bm{x}_{i})=\bm{0},
$$ 
$$
\frac{\nabla \bm{N}_i(\bm{x}_{i})}{\bm{N}_i(\bm{x}_{i})}=-\frac{\nabla P_{\bm{x}_\textrm{t}}(\bm{x}_{i})}{P_{\bm{x}_\textrm{t}}(\bm{x}_{i})}.
$$
\\
Therefore, since $||\nabla P_{\bm{x}_\textrm{t}}(\bm{x})/ P_{\bm{x}_\textrm{t}} (\bm{x})|| \leq J $,
$$
\bm{x}_i-\bm{x}_{i}^{\textrm{c}}= \bm{\Lambda}_i \cdot \frac{ \nabla P_{\bm{x}_\textrm{t}} (\bm{x}_{i})}{ P_{\bm{x}_\textrm{t}}(\bm{x}_{i})} = \mathcal{O}(1) \cdot \frac{1}{i^3} \frac{ \nabla P_{\bm{x}_\textrm{t}} (\bm{x}_{i})}{ P_{\bm{x}_\textrm{t}}(\bm{x}_{i})} = \mathcal{O}(1) \cdot \frac{1}{i^3} \cdot \bm{e}_i
$$
following the same idea as in Theorem \ref{thm1}, where $\bm{e}_i$ is the unit vector which have the same direction with 
$\nabla P_{\bm{x}_\textrm{t}}(\bm{x}_{i})/ P_{\bm{x}_\textrm{t}}(\bm{x}_{i})$.
Similarly, we can get 
$$
\bm{x}_{i+1}-\bm{x}_{i+1}^{\textrm{c}} = \mathcal{O}(1) \cdot \frac{1}{i^3} \cdot \bm{e}_{i+1}.
$$
Hence, 
$$
\begin{aligned}
\bm{x}_{i+1}&=\bm{x}_{i}-s_i(\bm{y}_i-\bm{y}_\textrm{t}) +\mathcal{O}(1) \cdot \frac{1}{i^3} \cdot \bm{p}_{i}\\
&=\bm{x}_{i}-s_i \nabla \bm{f}(\bm{z}_i)(\bm{x}_i-\bm{x}_\textrm{t}) + s_i \bm{\varepsilon}_i+\mathcal{O}(1) \cdot \frac{1}{i^3} \cdot \bm{p}_{i}
\end{aligned}
$$
with the mean value theorem, where $\bm{p}_i$ is a perturbation vector with a maximum norm of 2 and $\bm{z}_i$ lies between  $\bm{x}_i$ and $\bm{x}_\textrm{t}$.
Expanding the above equation, we can get
\begin{align}
& \bm{x_{i+1}}-\bm{x_\textrm{t}} \notag \\
=& \prod_{k \in \text{desc}(h, i)} \big(\bm{I} -s_k \nabla\bm{ f}(\bm{z_k})\big)(\bm{x_h}-\bm{x_\textrm{t}}) \label{eq34}\\
+&\sum_{k=h}^i\left(\prod_{u \in \text{desc}(k+1, i)}\big(\bm{I}-s_u \nabla\bm{ f}(\bm{z_u})\big)\right) \cdot s_k \bm{\varepsilon}_k \label{eq35}\\
+& \mathcal{O}(1) \cdot \sum_{k=h}^i\left(\prod_{u \in \text{desc}(k+1, i)}\big(\bm{I}-s_u \nabla\bm{ f}(\bm{z_u})\big)\right) \cdot \frac{1}{k^3} \cdot \bm{p}_{k} \label{eq36}
\end{align}
where $h$ is a large number s.t.\ $s_k<1 / M$, $\forall k \geq h$.
Since $m\bm{I} \prec \nabla\bm{ f} \prec M\bm{I}$ with $0 < m < M$,
$$
\left\|\prod_{k \in \text{desc}(h, i)} \big(\bm{I} -s_k \nabla\bm{ f}(\bm{z_k})\big)\right\| \leq \left\|\prod_{k \in \text{desc}(h, i)} \big(\bm{I} -s_k m\bm{ I}\big)\right\| \leq \prod_{k=h}^{i}(1- s_km ).
$$
Then, following this path, we can study the norm of the matrix products in Terms (\ref{eq34})--(\ref{eq36}) and derive that 
$$
\lim _{i \rightarrow \infty} \bm{x_{i+1}}-\bm{x_\textrm{t}} = \bm{0}
$$
following the idea of Theorem \ref{thm2}.
\end{proof}
\end{theorem}

\begin{theorem}
When Assumptions \ref{b1}--\ref{b6} hold, the multi-dimensional iterative Bayesian Robbins--Monro sequence has the same convergence rate as the standard multi-dimensional Robbins--Monro sequence in the long run.
\begin{proof}
We continue to investigate the norm of the product of matrices in Terms (\ref{eq34})--(\ref{eq36}). Then, following the same idea as in Lemma \ref{lem35} and Theorem \ref{thm36}, we can demonstrate that $\|E((\ref{eq35})^2)\| \gg \|(\ref{eq36})\|^2 $ when $i \to \infty$. 
Then, by sandwiching the norms of Terms (\ref{eq34})--(\ref{eq35}) as in Theorem \ref{thm36}, we can see that the multi-dimensional iterative Bayesian Robbins--Monro sequence has the same convergence rate as the standard multi-dimensional Robbins--Monro sequence in the long run.
\end{proof}
\end{theorem}

\section{TMS motor threshold estimation with IBRM} \label{tms}
\subsection{Background}

Transcranial magnetic stimulation (TMS) is a non-invasive neurostimulation technique that uses magnetic pulses to write signals into specific neural circuits of the brain \cite{bhattacharya2022overview}. By placing a coil over the scalp, a focused magnetic field pulse induces electrical currents in the underlying cortical tissue and triggers neural signals in response \cite{goetz2017development}. TMS has become a valuable tool for research in experimental brain research and for clinical treatment, particularly in the therapy of brain disorders, such as clinical depression \cite{sonmez2019accelerated}.   
All procedures involve an initial dosage individualisation through the titration of the motor threshold.
The motor threshold is defined as the minimum intensity of stimulation to elicit a motor evoked potential (MEP) with a specific medial amplitude (usually $\geq 50~$µV peak-to-peak) in a target muscle (typically in the hand) \cite{westin2014determination}.The response is highly variable but the expectation is monotonically increasing with the stimulation strength \cite{ma2024extraction}. The individual's motor threshold establishes the safety limits as large stimulation strength in various procedures can trigger seizures. Furthermore, the threshold individualises all subsequent stimulation strengths for treatment \cite{ma2023correlating}.   

The Robbins--Monro sequence serves for detecting the motor threshold with few stimuli \cite{liu2025prior}.
The function $f$ describes the relationship between the stimulation strength (input $x$) and the response.

The field uses two fundamental variants. The first uses the continuous (\textit{analog}, ACS) log-transformed response (output $y$) and a target level $y_{\textrm{t}}$ of ln(50\,µV) to find $x_{\textrm{t}}$ as the threshold.
The binary (\textit{digital}, DCS) version only differentiates between response ($y_i-y_\textrm{t}\geq 0$), which leads to a step of $s_i$ down ($x_{i+1}=x_i-s_i$), or no response ($y_i-y_\textrm{t}< 0$), which leads to a step of $s_i$ up ($x_{i+1}=x_i+s_i$).


We compare the standard Robbins--Monro sequence (DCS and ACS), prior-information Robbins--Monro sequence (see Eq.\! (\ref{def})) of both variants (DCS-PI and ACS-PI), and iterative Bayesian Robbins--Monro sequence (DCSu and ACSu). Note that $s_i=s/i$ for DCS, ACS, DCS-PI and ACS-PI, and $s_i=s$ for DCSu and ACSu.

\subsection{Result and discussion}

\begin{figure}[H]
    \centering
    \includegraphics[scale=0.38]{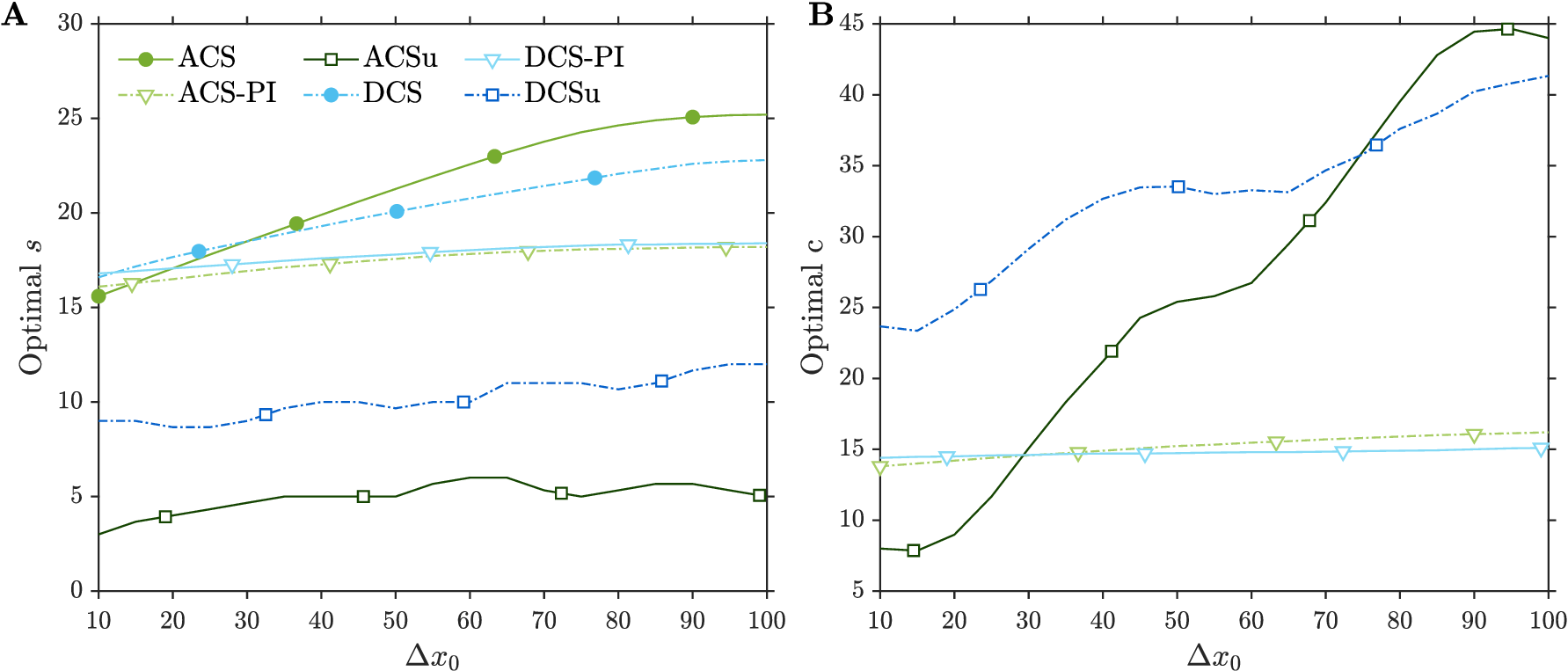}
    \caption{\textbf{A}. The relationship between the optimal $s$ and $\Delta x_0=x_0-x_\textrm{t}$ for all six sequences. 
    \textbf{B}. The relationship between the optimal $c$ and $\Delta x_0=x_0-x_\textrm{t}$ for the four sequences with prior information. 
    The results in subplots A and B are averaged over 100,000 $s$ and $c$ values, from 600,000 randomly selected $s$ and $c$ that achieve the smallest absolute error at step 30 ($|\Delta x_{30}| = |x_{30}-x_\textrm{t}|$). }
    \label{f1}
\end{figure}

We initially identify the optimal individual $s$ and $c$ values for each method.
We estimate the TMS threshold using a sample of 25,000 virtual subjects and each subject is repeated 24 times. The subject data are generated based on a detailed stochastic model of brain stimulation \cite{goetz2019statistical}. The prior distribution of TMS thresholds ($x_\textrm{t}$) among these subjects follows approximately a normal distribution $\mathcal{N}(65,10^2)$, with an average threshold at $x=65$ machine output and a standard deviation of $10$ \cite{goetz2019statistical,wang2023three}.
The values of $s$ and $c$ are varied randomly during these runs of each algorithm. The optimum values for each algorithm are chosen as the values with the lowest error after 30 stimuli, i.e., evaluations of function $f$ (Figure \ref{f1}).

\begin{figure}[H]
    \centering
    \hspace*{-0.8cm}
    \includegraphics[scale=0.33]{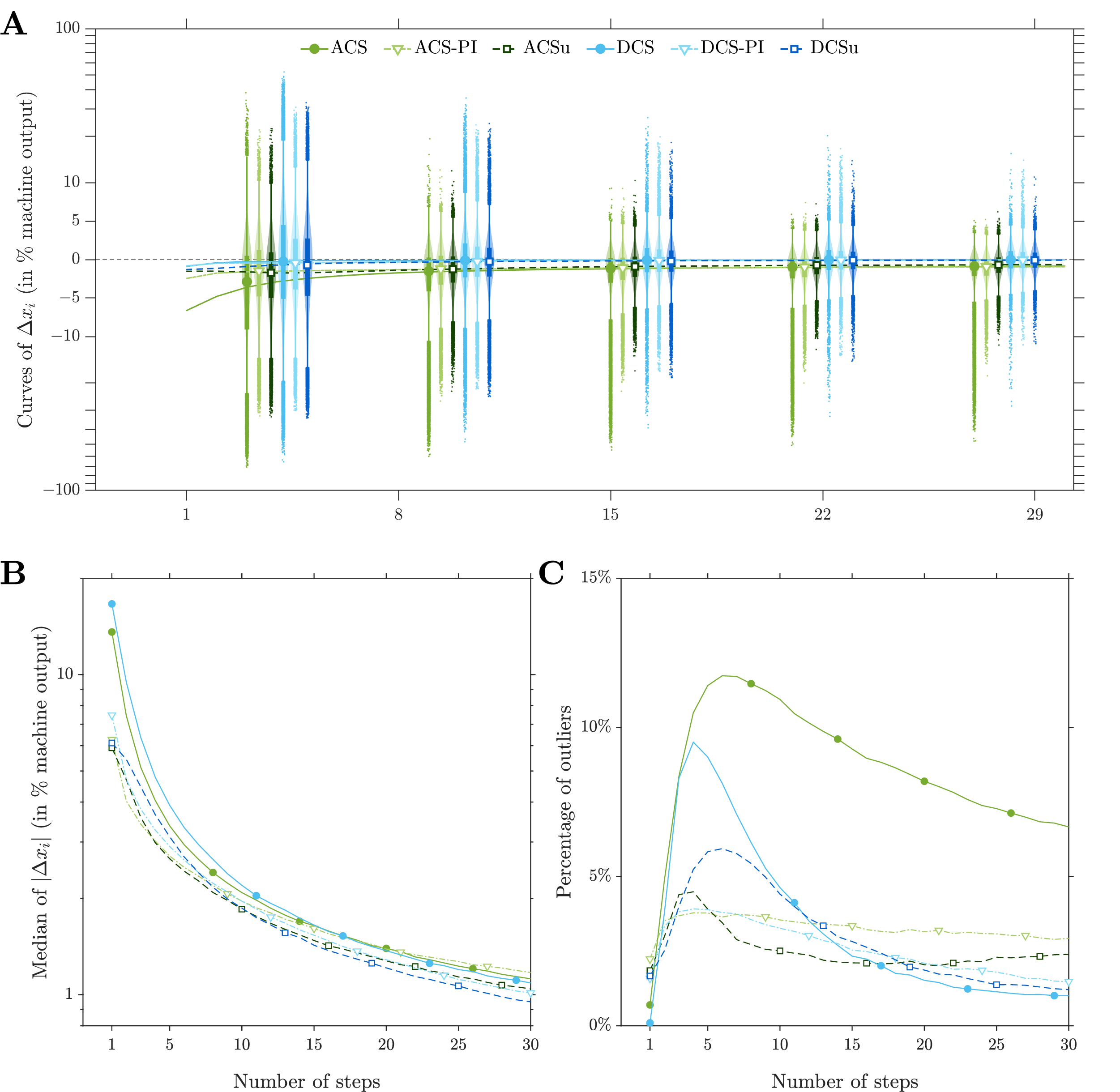}
    \caption{Convergence of the six sequences for 25,000 virtual subjects. 
    (\textbf{A}) Statistics of threshold estimation errors ($\Delta x_i = x_i-x_\textrm{t}$) over 30 steps (mixed log--lin scale, lines with markers: median error, bars: quartiles, whiskers: upper and lower adjacent values that are 1.5-times the interquartile range below and above the 25-th and 75-th percentiles, respectively, and dots: outliers outside the upper and lower adjacent values).
    (\textbf{B}) Median absolute error $|\Delta x_i|$.
    (\textbf{C}) Percentage of outliers.}
    \label{f3}
\end{figure}

In the subsequent test runs, we set the initial step size $s$ to 20 for ACS and DCS, 17 for ACS-PI and DCS-PI, 5 for DCSu, and 10 for ACSu. The parameter $c$ was set to 15 for ACS-PI as well as DCS-PI, and 30 for ACSu as well as DCSu. Figure \ref{f3} displays the results.

Both DCSu and ACSu excel in speed and robustness compared to other approaches.
DCSu and ACSu reduce the spread and percentage of outliers compared to the corresponding standard Robbins--Monro methods (DCS and ACS), and ACSu shows a much narrower spread of outliers than ACS-PI in the long run (see Figure \ref{f3} A \& C). 

Moreover, DCSu and ACSu achieve consistently lower median absolute deviations $|\Delta x_i|$ compared to other methods (see Figure \ref{f3} B).

Thus, the here suggested Bayesian Robbins--Monro variants promise faster threshold measurement with fewer stimuli and a lower risk of outliers, which otherwise could increase the seizure risk in subsequent procedures.

\section{Conclusion} \label{con}

This paper presented an iterative Bayesian Robbins--Monro sequence which adaptively and iteratively updates a prior distribution with new observations to enhance the convergence speed and stability of standard Robbins--Monro. 
Theoretical analysis confirmed its almost sure convergence and established its rate. A practical implementation for brain stimulation threshold estimation demonstrated its robustness, accuracy, and speed. The new method reached the accuracy level of the standard Robbins--Monro at 30 steps already after 22 steps. Furthermore, it had fewer than half of the outliers. In brain stimulation, the number of function evaluations, i.e., stimuli, needs to be as low as possible to minimise side-effects, while the safety of the subsequent procedures depends on the detection accuracy of the threshold level.

\section*{Competing interests}
No competing interest is declared. No external funds, grants, or other support were received. The authors have no relevant financial or non-financial interests to disclose.

\section*{Author contributions statement}

S.M.G. conceived, supervised, secured funding for the study, and provided the sequence for introducing prior information into the RM iteration. S.M.G. also provided part of the code. S.M.G. and S.L. conceived and sketched concepts for the proofs. S.L. completed all the proofs in this paper. S.L. further studied and characterized the sequence and its performance. S.M.G. and K.M. checked the proofs. S.L. and K.M. also wrote the rest of the code for numerical analysis and visualization. S.L. and S.M.G. wrote the manuscript, K.M. edited the text and figures, and all authors reviewed, commented on, and approved the final version of the manuscript.

\section*{Acknowledgments}
The authors thank the anonymous reviewers for their valuable suggestions. 

\section*{Data availability}
The data will be made available upon request.


\bibliography{sn-article}

\end{document}